\documentclass[11pt]{article}

\usepackage{tikz}
\usetikzlibrary{arrows,decorations.pathmorphing,
                backgrounds,positioning,fit,petri}
\usepackage{amsmath}
\usepackage{amssymb}
\usepackage{amsbsy}
\usepackage{amsfonts}
\usepackage{mathtools}
\usepackage{extarrows}
\usepackage{bbm}
\usepackage{comment}
\usepackage{xcolor}
\newtheorem{theorem}{Theorem}[section]%
\newtheorem{lemma}[theorem]{Lemma}%
\newtheorem{prop}[theorem]{Proposition}%

\newenvironment{pf}{\medskip\noindent{Proof:}
  \hspace{-.5cm}      \enspace}{\hfill \qed \newline \smallskip}

   \def\z{\zeta}

\def\det{{\rm det}\,}

\def\f{\noindent}

\def\mod{\hbox{\rm mod}\;}

\newcommand{\qed}{\mbox{\raisebox{0.7ex}{\fbox{}}} \vspace{4truemm}}

\begin{document}

\begin{center}
{\Large{\textbf{On complexity of cyclic coverings of graphs}}}
\end{center}

\vskip 5mm {\small
\begin{center}
{\textbf{Y.~S.~Kwon,}}\footnote{{\small\em Department of
Mathematics, Yeungnam University, Korea}}
{\textbf{A.~D.~Mednykh,}}\footnote{{\small\em Sobolev Institute of
Mathematics, Novosibirsk State University, Russia }} {\textbf{I.~A.~Mednykh,}}\footnote{{\small\em Sobolev
Institute of Mathematics, Novosibirsk State University, Russia }}
 \end{center}}

\title{ \vspace{-1.2cm}
On complexity of cyclic coverings of graphs\thanks{Supported by }}

\begin{abstract}
By complexity of a finite graph we mean the number of spanning trees in the graph. The aim of the present paper is to give a new approach  for counting complexity $\tau(n)$ of cyclic $n$-fold coverings of a graph. We give an explicit analytic formula for $\tau(n)$ in terms of Chebyshev polynomials and find its asymptotic behavior as $n\to\infty$  through the Mahler measure of the associated voltage polynomial. We also prove that
$F(x)=\sum\limits_{n=1}^\infty\tau(n)x^n$ is a rational function with integer coefficients.
\bigskip

\f\textbf{AMS classification:} 05C30, 39A10

\end{abstract}

\section{Introduction}

A {\it spanning tree} $T$ in a graph $G$ is a subgraph that
is a tree containing all   the vertices of $G$. The number $\tau(G)$ of spanning trees in a connected graph $G$ (or complexity of $G$) is a well studied invariant. In the simplest cases it can be calculated explicitly. For instance, if $G=C_n$ is the cycle graph on $n$ vertices, then $\tau(G)=n,$ if $G=K_n$ is a complete graph   on $n$ vertices, then  $\tau(G)=n^{n-2}$ (the Caley's formula), if $G=Q_n$ is the $n$-dimensional hypercube, then $\tau(G)=2^{2^n-n-1}\cdot\prod_{k=2}^{n}k^{{n}\choose{k}}$ (\cite{HaHaWu}).
More complicated formulas for the number of spanning trees are known for some special graphs,
such as the wheel, fan, ladder, M\"obius ladder
\cite{BoePro}, grids \cite{NP04}, lattices \cite{SW00}, prism and
anti-prism \cite{SWZ16} and some self-similar graphs  (\cite{CCY07}, \cite{AD11}). In many important families of graphs,  the complexity of graphs in the families is expressed is terms of the Chebyshev
polynomials, while its asymptotical behavior is defined by the Mahler measure of the associated polynomial. This takes place for the circulant graphs \cite{ZhangYongGol},
\cite{ZhangYongGolin}, \cite{XiebinLinZhang}, \cite{MedMed17}, \cite{MedMed18}, the generalized Petersen graphs \cite{KwonMedMed}, the $I$-graphs \cite{Ilya}, and some other families of graphs \cite{Louis15}, \cite{SW16}, \cite{SW17}, \cite{AbrBaiMed}.

In this paper, we consider an infinite family $H_n$ of $n$-fold cyclic coverings of a finite connected graph $H$ defined by a voltage assignment on $H$. We provide an explicit formula for the number of spanning trees $\tau(n)$ in $H_n$ in terms of Chebyshev polynomials  and find its asymptotic behavior as $n\to\infty$ through the Mahler measure  of the associated voltage polynomial.  Also we prove that
$F(x)=\sum\limits_{n=1}^\infty\tau(n)x^n$  is a rational function with integer coefficients.

\section{Basic definitions and preliminary facts}\label{basic}

Consider a connected finite graph $G,$ allowed to have multiple edges and loops.   We denote the vertex and edge set of $G$ by $V(G)$ and $E(G),$ respectively. Given $u, v\in V(G),$ we set $a_{uv}$ to be equal to the number of edges between vertices $u$ and $v$ if $u \neq v$; $2$ times the number of loops incident to $u$ if $u=v$.
The matrix $A=A(G)=\{a_{uv}\}_{u, v\in V(G)}$ is called \textit{the adjacency matrix} of the graph $G.$ The degree $d_v$ of a vertex $v \in V(G)$ is defined by $d_v=\sum_{u\in V(G)}a_{uv}.$ Let $D=D(G)$ be the diagonal matrix indexed by the elements of $V(G)$ with $d_{vv} = d_v.$ The matrix $L=L(G)=D(G)-A(G)$ is called \textit{the Laplacian matrix}, or simply \textit{Laplacian}, of the graph $G.$ In what follows, by $I_n$ we denote the identity matrix of order $n.$

We call an $n\times n$ matrix {\it circulant,} and denote it by $circ(a_0, a_1,\ldots,a_{n-1})$ if it is of the form
$$circ(a_0, a_1,\ldots, a_{n-1})=
\left(\begin{array}{ccccc}
a_0 & a_1 & a_2 & \ldots & a_{n-1} \\
a_{n-1} & a_0 & a_1 & \ldots & a_{n-2} \\
  & \vdots &   & \ddots & \vdots \\
a_1 & a_2 & a_3 & \ldots & a_0\\
\end{array}\right).$$

Recall \cite{PJDav} that the eigenvalues of matrix $C=circ(a_0,a_1,\ldots,a_{n-1})$ are given by the following simple formulas $\lambda_j=p(\varepsilon^j_n),\,j=0,1,\ldots,n-1$
where $p(x)=a_0+a_1 x+\ldots+a_{n-1}x^{n-1}$ and $\varepsilon_n$ is an order $n$ primitive root of the unity.
Moreover, the circulant matrix $C=p(T),$ where $T_{n}=circ(0,1,0,\ldots,0)$ is the matrix representation of
the shift operator $T_{n}:(x_0,x_1,\ldots,x_{n-2},x_{n-1})\rightarrow(x_1, x_2,\ldots,x_{n-1},x_0).$
For any $i=0,\ldots, n-1$, let ${{\bf v}_i} =(1,\varepsilon_n^i ,\varepsilon_n^{2i},\ldots,\varepsilon_n^{(n-1)i})^t$ be a column vector of length $n.$
We note that all $n\times n$ circulant matrices share the same set of linearly independent eigenvectors $\textbf{v}_0, \textbf{v}_1, \ldots, \textbf{v}_{n-1}.$
Hence, any set of $n\times n$ circulant matrices can be simultaneously diagonalizable.


For a subgraph $H$ of $G$, if each connected component of $H$ is a cycle,
loop or an edge, we call $H$ \emph{basic subgraph} of $G$. Denote the set
of all basic subgraphs of $G$ by ${\cal B}(G)$. For any basic subgraph $B$
of $G$, let $B_1$ be the subgraph of $B$ composed of all components of $B$
which is an edge.  By ${\cal C}(G)$ and ${\cal M}(G)$, we denote the set
of all loops or cycles of $G$ and the set of all matchings of $G$, respectively.
Note that the empty subgraph of $G$ also belongs to ${\cal M}(G)$. For any
$k=1,\ldots,\ell(G)$, let ${\cal C}_k(G)$ be the set of basic subgraphs $B$
such that the number of components of $B$ is $k$ and any component of $B$
is a loop or a cycle, where $\ell(G)$ is the maximum number of components
of basic subgraphs of $G$ whose each component is a loop or a cycle. Note
that ${\cal C}(G)={\cal C}_1(G)$. For a loop or a cycle $C$, let ${\cal C}_k(G,C)$
be the set of basic subgraphs in ${\cal C}_k(G)$ containing $C$ as a component.

\section{Cyclic coverings and their Laplacians}
In this section, we describe the cyclic coverings of a finite connected graph  and investigate the structure of their Laplacian.

Let $H$  be a finite connected graph with  the vertices
$v_1,v_2,\ldots,v_r,$  possibly with loops and multiple edges.
Denote by $a_{i,j}$ the number of edges between $v_i$ to $v_j$ if
$i\neq j$; $2$ times the number of loops incident to $v_i$ if $i=j$.
 Then the adjacency matrix of the graph $H$ is
$$A(H)=
\left(\begin{array}{cccc}
a_{1,1} & a_{1,2}   & \ldots & a_{1,r} \\
a_{2,1} & a_{2,2}   & \ldots & a_{2,r} \\
\vdots   & \vdots    & \ddots & \vdots \\
a_{r,1} & a_{r,2}   & \ldots & a_{r,r}\\
\end{array}\right).$$
To  describe the general construction of a cyclic covering of $H,$ we consider the
cyclic group  $\Gamma$ (finite or infinite)  and  we suppose that all edges
of $H,$ including loops are oriented in  two possible ways and
denote the corresponding digraph by $\overrightarrow{H}$. For any
oriented edge $e$ of  $\overrightarrow{H}$, we denote the oppositely
oriented edge coming from the same edge with $e$ by $e^{-1}$.

By making use of the voltage technique developed in
\cite{GrossTucker}, we assign to each oriented edge $e$ of
$\overrightarrow{H}$ a {\it voltage} $\alpha(e)$, namely    some
element of the group $\Gamma$. We follow the rule that the sum of
voltages assigned to two oriented edges induced by the same edge is
zero, namely $\alpha(e) + \alpha(e^{-1})=0$. For any $i, j=1,
\ldots, r$  and for any $k=1, \ldots, a_{i,j}$,  we denote the
voltage of the $k$-th oriented edge from $v_i$ to $v_j$ by
$\alpha_k(i,j)$.
 For different $i$ and $j$, we assume that  $\alpha_k(j,i) =-\alpha_k(i,j)$. Since each loop  has two opposite orientations,  one can take
the voltage $\alpha_k(i,i)$ for one of them (arbitrary chosen) and
$-\alpha_k(i,i)$  for another.   For any cycle
$C$ of $H$, let $C^{+}$ be an orientation of $C$ such that it is a
directed cycle. Let $C^{-}$ opposite directed cycle with $C^{+}.$ By
$\alpha(C^{+}),$ we denote the sum of all voltages assigned to
directed edges belonging to $C^+.$  Note that
$\alpha(C^{-})=-\alpha(C^{+}).$ Fix a spanning tree $T$ in $H$. For
any $e \notin E(T)$, $T+e$ contains a unique cycle. Denote the cycle
by $C_e.$ Now it is well known that the covering of $H$ derived by the voltage assignment $\alpha$ is
connected if and only if the voltages $\{\alpha(C_e^+)\ | \ e \notin
E(T)\}$ generate the full group $\Gamma$. For our convenience, we assume that all voltages assigned to
oriented edge in $\overrightarrow{H}$ are integers and  the voltages
$\{\alpha(C_e^+)\ | \ e \notin E(T)\}$ generate the full group
$\mathbb{Z}$, namely,   there exist cycles $C_1,  \ldots, C_s$
 such that $k_1\alpha(C_1^+)+ \cdots +k_s\alpha(C_s^+) = 1$ for some integers
$k_1,\ldots, k_s$.

 Given a positive integer $n$,  we consider each voltage modulo $n$. Denote by $H_n = H_{n}^{\alpha}$ the $n$-fold cyclic covering of $H$  derived by  the voltage assignment $\alpha=\{\alpha_k(i,j),\,i, j=1,
\ldots, r,\,k=1, \ldots, a_{i,j}\}.$  According to the voltage theory \cite{GrossTucker}, the graph
$H_{n}^{\alpha}$ has the  set of vertices
$u_{i,\,s},\,i=1,2,\ldots,r,\,s=1,2,\ldots,n,$  and the  set of
edges $v_{i,\,s}v_{j,\,s+\alpha_k(i,j)},$  where $\,i,j=1,2,\ldots,
r,\, k = 1,2,...,a_{i,j},$ and the second indexes are taken  $\mod
n.$

We emphasize  that all cyclic $n$-fold coverings of the graph $H$  can be  obtained in such a way. In particular, this construction covers
the circulant graphs, Haar graphs, $I$-graphs, $Y$-graphs, $H$-graphs and many other
families of graphs.

The   adjacency matrix $A(H_{n})$ for the graph $H_{n}$ is a block matrix
$$
A(H_{n})=\left(\begin{array}{cccc}
B_{11} & B_{12}  & \ldots & B_{1r} \\
B_{21} & B_{22}  & \ldots & B_{2r} \\
  \vdots &   & \ddots & \vdots \\
B_{r1} & B_{r2}  & \ldots & B_{rr}\\
\end{array}\right),$$  where
the $(i,j)$-th block entry $B_{ij}=B_{ij}(n),\,i,j=1,2,\ldots,r$ is given by the  formulas \begin{equation}\label{block}B_{ij}=\sum\limits_{k=1}^{a_{i,j}}T_{n}^{\alpha_k(i,j)},\,i\neq j,   \,\,\, and\,\,\,B_{ii}=\sum\limits_{k=1}^{a_{i,i}}(T_{n}^{\alpha_k(i,i)}+ T_{n}^{-\alpha_k(i,i)}).\end{equation}

Here, $B_{ij}$  is an  $ n\times  n$ matrix whose
$(s,s^{\prime})$-th entry   is equal to the number of edges between
vertices $v_{i,\,s}$  and $v_{j,\,s^{\prime}}.$
 We also represent the degree matrix  $D(H_{n})$  of the graph $H_{n}$  in the form
 $$
D(H_{n})=\left(\begin{array}{cccc}
D_{11} & 0  & \ldots & 0 \\
0 & D_{22}  & \ldots & 0 \\
  \vdots &   & \ddots & \vdots \\
0 & 0  & \ldots & D_{rr}\\
\end{array}\right),$$
where  $D_{ii}$ is a diagonal  $n\times n$  matrix whose $s$-th
entry $d_{i,s}$  is the valency of vertex $v_{i,s}$. Since the
vertex $v_{i,s}$  is connected by an oriented edge with vertices $v_{j,\,s+\alpha_k(i,j)},$ {} $\,j=1,2,\ldots, r,\, k =
1,2,...,a_{i,j},\,i\neq j$  and by $2a_{i,i}$ oriented edges with
a fibre of $v_i$, we have
\begin{equation}\label{degree}d_{i,s}=a_{i,i}+\sum\limits_{j=1}^ra_{i,j},\end{equation}
which equals to the degree $d(v_i)$ of $v_i$. By (\ref{degree}), the
valency $d_{i,s}$ does not depend on $s,$ so
\begin{equation}\label{matrixdegree}D_{ii}= d(v_i) I_n.\end{equation}
Recall that the Laplacian matrix  $L(H_{n})=D(H_{n})-A(H_{n}).$
Hence, the $(i,j)$-th block entry of the  matrix
$L(H_{n})=(L_{ij})_{i,j=1,2,\ldots,r}$   can be written as
\begin{equation}\label{laplacian}
L_{ij}=\begin{cases} d(v_i) I_n-B_{ii}
& \text{ if } i =j \\
-B_{ij} & \text{ if } i \neq j .
\end{cases} \end{equation}

\section{Kirchhoff theorem and the number of spanning trees in $H_{n}$}

By the classical Kirchhoff theorem, the  number of spanning trees of
the graph $H_{n}$ is equal  to the product of non-zero eigenvalues
of $L(H_{n})$ divided by the number of vertices $r\,n$  of $H_{n}.$
 To find spectrum of $L(H_{n})$  we have to solve the linear system of equations
 \begin{equation}\label{linearsystem}
 \left(\begin{array}{cccc}
L_{11}-\lambda I_n & L_{12}  & \ldots & L_{1r} \\
L_{21} & L_{22}-\lambda I_n  & \ldots & L_{2r} \\
  \vdots &   & \ddots & \vdots \\
L_{r1} & L_{r2}  & \ldots & L_{rr}-\lambda I_n\\
\end{array}\right)\left(\begin{array}{c}X_1\\X_2\\\vdots\\X_r\end{array}\right)=\left(\begin{array}{c}0\\ 0\\\vdots\\ 0\end{array}\right),\end{equation}
where  $X_i,\,i=1,2,\ldots,r$  is an $n$-dimensional column vector.

We associate   with   $L(H_{n})$ the following   Laurent polynomial

\begin{equation}\label{polynomial}
P(z,\lambda)=\det\left(\begin{array}{cccc}
p_{11}(z)-\lambda & p_{12}(z)  & \ldots & p_{1r}(z) \\
p_{21}(z) & p_{22}(z)-\lambda  & \ldots & p_{2r}(z) \\
  \vdots &   & \ddots & \vdots \\
p_{r1}(z) & p_{r2}(z)  & \ldots & p_{rr}(z)-\lambda\\
\end{array}\right),\end{equation} where $p_{ii}(z)=d(v_i)-\sum\limits_{k=1}^{a_{i,i}}(z^{\alpha_{k}(i,i)}  + z^{-\alpha_{k}(i,i)} )$  and
$p_{ij}(z)=-\sum\limits_{k=1}^{a_{i, j}}z^{\alpha_k(i,j)}$ for $i
\neq  j.$ One can easily check that $L_{ij}=p_{ij}(T_n),$   where $T_n={\rm circ}\,(0,1,0,\ldots,0)$ is an $n\times n$ circulant matrix. We  also consider the polynomial

\begin{equation}\label{Vpolynomial}
P(z)=\det\left(\begin{array}{cccc}
p_{11}(z)  & p_{12}(z)  & \ldots & p_{1r}(z) \\
p_{21}(z) & p_{22}(z)   & \ldots & p_{2r}(z) \\
  \vdots &   & \ddots & \vdots \\
p_{r1}(z) & p_{r2}(z)  & \ldots & p_{rr}(z)\\
\end{array}\right). \end{equation}\bigskip  We will refer to $P(z)$ as a {\it voltage polynomial}  of the cyclic covering $H_n.$ We have $P(z,0)=P(z).$

To solve the system~(\ref{linearsystem}) we make the following
observation.  The eigenvalues of circulant matrix $T_{n}$ are
$\varepsilon_{n}^{j},\,j=0,1,\ldots,n-1,$ where
$\varepsilon_n=e^\frac{2\pi i}{n}.$ Since all of them are distinct,
the matrix $T_{n}$ is conjugate to the diagonal matrix
$\mathbb{T}_{n}=diag(1,\varepsilon_{n},\ldots,\varepsilon_{n}^{n-1})$
with diagonal entries
$1,\varepsilon_{n},\ldots,\varepsilon_{n}^{n-1}$.

To find spectrum of $L(H_{n}),$ without loss of generality, one can
assume  that $T_{n}=\mathbb{T}_{n}.$ Then all blocks of $L(H_{n})$
are diagonal matrices. Hence, the linear system~(\ref{linearsystem})
splits into $n$ independent linear systems

\begin{equation}\label{splitlinearsystem}
 \left(\begin{array}{cccc}
L_{11}^s-\lambda & L_{12}^s  & \ldots & L_{1r}^s\\
L_{21}^s & L_{22}^s-\lambda  & \ldots & L_{2r}^s \\
  \vdots &   & \ddots & \vdots \\
L_{r1}^s & L_{r2}^s  & \ldots & L_{rr}^s-\lambda\\
\end{array}\right)\left(\begin{array}{c}X_{1,s}\\X_{2,s}\\\vdots\\X_{r,s}\end{array}\right)=\left(\begin{array}{c}0\\ 0\\\vdots\\ 0\end{array}\right),\,s=1,2,\ldots,n,\end{equation} where $X_{i,s}$ is the $s$-th component of the vector $X_i$ and $L_{i,j}^s,$ is the $s$-th entry of the diagonal matrix $L_{ij}.$ Indeed, we have $L_{ij}^s=p_{ij}(\varepsilon_n^s),$  where $p_{ij}(z)$  are the same as in formula (\ref{polynomial}).
Hence, to find eigenvalues we have to solve  $n$ algebraic equations  of degree $r$

\begin{equation}\label{algebraic}
P(\varepsilon_n^s,\lambda)=0,\, s=1,2,\ldots,n.
\end{equation}
Given $s,$ we denote the solutions of the $s$-th equation by $\lambda_{s,i},\,i=1,2,\ldots,r.$  By the Vietta theorem  we get
\begin{equation}\label{vietta}
\lambda_{s,1}\lambda_{s,2}\ldots\lambda_{s,r}=P(\varepsilon_n^s,0)=P(\varepsilon_n^s).
\end{equation}
 In particular case  $s=n,$  the  linear system~(\ref{splitlinearsystem})  gives the eigenvalues and eigenvectors of the graph $H.$  Since graph $H$  is connected, it has only one zero eigenvalue, say $\lambda_{n,1}.$  By the Kirchhoff theorem, the product of the  non-zero eigenvalues of $H$  is equal to the number of vertices of $H$  multiplied by its complexity.  Thus \begin{equation}\label{kirhhoffH}
\prod_{i=2}^r\lambda_{n,i}=r\,\tau(H).
\end{equation}

Now, we use the Kirchhoff theorem to find complexity of graph $H_{n}.$  Again, since $H_{n}$  is connected it has exactly one zero eigenvalue $\lambda_{n,1}=0.$  Taking into account~(\ref{kirhhoffH})  we get

\begin{equation}\label{kirhhoffHmn}
\tau(H_{n})=\frac{1}{r\,n}\prod_{i=2}^r\lambda_{n,i}\cdot\prod\limits_{s=1}^{n-1}\lambda_{s,1}\lambda_{s,2}\ldots\lambda_{s,r} =\frac{\tau(H)}{n}\prod\limits_{s=1}^{n-1}\lambda_{s,1}\lambda_{s,2}\ldots\lambda_{s,r}.
\end{equation}

Substituting (\ref{vietta}) into (\ref{kirhhoffHmn}), we obtain

\begin{equation}\label{finalkirhhoffHmn}
\tau(H_{n})=\frac{\tau(H)}{n}\prod\limits_{s=1}^{n-1}P(\varepsilon_n^s).
\end{equation}

We get the following preliminary result.
\bigskip

\begin{prop}\label{theorem0}
The number of spanning trees in the graph $H_{n}$ is given by the formula~(\ref{finalkirhhoffHmn}), where $P(z)$  is the voltage polynomial of $H_{n},$ $\tau(H)$ is the number of spanning trees in the graph $H$ and $\varepsilon_n$ is the $n$-th primitive root of the unity.
\end{prop}

\section{Voltage polynomial and its properties}

Now we are going to investigate the basic properties of the voltage polynomial $P(z).$
\bigskip

Note that $P(1) = \det L(H)=0$ and

\begin{eqnarray} \label{Sach}
P(z) &=& \sum_{\sigma \in S_r} sgn(\sigma)p_{1,\sigma(1)}
(z)p_{2,\sigma(2)}(z)\cdots p_{r,\sigma(r)}(z) \\
&=& \label{Sachh}\sum_{B \in {\cal B}(H)}(-1)^{|E(B_1)|} \left(\prod_{C \in {\cal
C}(B)} (-z^{\alpha(C^+)}-z^{\alpha(C^-)})\right) \left(\prod_{v \in
V(H)\setminus V(B)}d(v) \right)
\end{eqnarray}

For our convenience, we use $p_{i,j}$  instead of $p_{i,j}(z).$  Note that if $\sigma(i)=i$, then
$p_{i,\sigma(i)}=p_{i,i}=d(v_i)-\sum\limits_{k=1}^{a_{i,i}}(z^{\alpha_{k}(i,i)}
+ z^{-\alpha_{k}(i,i)} )$. In the polynomial expansion ~(\ref{Sach}) of $P(z)$, if
we choose $d(v_i)$ in $p_{i,i}$, then it contributes to the term $\left(\prod_{v
\in V(H)\setminus V(B)}d(v) \right)$; and if we choose
$-z^{\alpha_{k}(i,i)}$ or $-z^{-\alpha_{k}(i,i)}$ in $p_{i,i}$, then it
contributes to the polynomial $\left(\prod_{C \in {\cal C}(B)}
(-z^{\alpha(C^+)}-z^{\alpha(C^-)})\right)$.

Suppose that $\sigma(i) =j$ and $\sigma(j)=i$ for some $i,j$ with $i
\neq j$ and there is an edge between the vertices $v_i$ and $v_j.$  The transposition $(i ,j)$ contributes $-1$ to $sgn(\sigma)$
and
$$p_{i, \sigma(i)}p_{j, \sigma(j)} = p_{i, j}p_{j, i}= \left(-\sum\limits_{k_1 =1}^{a_{i,
j}}z^{\alpha_{k_1}(i,j)}\right) \times
\left(-\sum\limits_{k_2=1}^{a_{i, j}}z^{\alpha_{k_2}(j,i)}\right).$$
If there is no edge between $v_i$ and $v_j$, then $p_{i,
\sigma(i)}p_{j, \sigma(j)} =0$. When there is at least one edge
between $v_i$ and $v_j$, in the polynomial expansion of $P(z)$, if
we choose $-z^{\alpha_{k_1}(i,j)}$ in $p_{i, j}$ and $-z^{\alpha_{k_2}(j,i)}$ in $p_{j, i}$ with
$k_1 = k_2$, then it contribute to the term $(-1)^{|E(B_1)|}$ with
the effect of $sgn(\sigma)$; and if we choose
$-z^{\alpha_{k_1}(i,j)}$ in $p_{i, j}$ and $-z^{\alpha_{k_2}(j,i)}$  in $p_{j, i}$  with $k_1 \neq
k_2$, then it contribute to the polynomial $\left(\prod_{C \in {\cal
C}(B)} (-z^{\alpha(C^+)}-z^{\alpha(C^-)})\right)$ with the effect of
$sgn(\sigma)$.

Let the cyclic permutation $(i_1, i_2, \ldots, i_t)$ be a cycle of length $t\ge3$ in
the decomposition of $\sigma$ into disjoint cycles. If there is no
cycle containing $v_{i_1}, v_{i_2}, \ldots, v_{i_t}, v_{i_1}$ in
this order, then $p_{i_1,i_2}p_{i_2, i_3}\cdots
p_{i_{t-1},i_t}p_{i_t, i_1}=0$. If there is at least one such cycle,
then in the polynomial expansion of $P(z)$, $p_{i_1,i_2}p_{i_2,
i_3}\cdots p_{i_{t-1},i_t}p_{i_t, i_1}$ contributes  to the
polynomial $\left(\prod_{C \in {\cal C}(B)}
(-z^{\alpha(C^+)}-z^{\alpha(C^-)})\right)$ with the effect of
$sgn(\sigma)$. So we have the equation ~(\ref{Sachh}).

Note that the constant term of $P(z)$ is ${\displaystyle \sum_{M \in
{\cal M}(H)}(-1)^{|E(M)|}\prod_{v \in V(H)\setminus V(M)}d(v)}$ and
for any $k=1,\ldots,\ell(H)$ and for any $B \in {\cal C}_k(H)$, the
coefficient of $z^{\sum_{C \in {\cal C}(B)}  \alpha(C^{\pm})}$ is
$$(-1)^k\sum_{M \in {\cal
M}(H-B)}(-1)^{|E(M)|}\prod_{v \in V(H-B-M)}d(v),$$ where
$\alpha(C^{\pm})$ means $\alpha(C^{+})$ or $\alpha(C^{-})$. Denote
this coefficient by $c_{B}$ and let $c_0$ be the constant term of
$P(z)$. Now we have

\begin{eqnarray}\label{newform}
P(z) &=&  \sum_{M \in {\cal M}(H)}(-1)^{|E(M)|}\prod_{v \in
V(H)\setminus
V(M)}d(v) \nonumber  \\
&+&\sum_{k=1}^{\ell(H)}(-1)^{k}\sum_{B \in {\cal
C}_k(H)}\left(\prod_{C \in {\cal C}(B)}(z^{\alpha(C^+)}
+z^{\alpha(C^-)})\right)\left(\sum_{M \in {\cal
M}(H-B)}(-1)^{|E(M)|}\prod_{v \in V(H-B-M)}d(v) \right)\nonumber  \\
&=& c_0 +\sum_{k=1}^{\ell(H)}\sum_{B \in {\cal
C}_k(H)}c_B\left(\prod_{C \in {\cal C}(B)}(z^{\alpha(C^+)}
+z^{\alpha(C^-)})\right).
\end{eqnarray}

 For later use, we need
the following lemma.

\begin{lemma}\label{positive-det}
Let $a_{i,j},(a_{i,i}=0),\,i,j=1,2,\ldots,m$ be non-negative numbers. Let $$D(x_1,x_2,\ldots,x_m)=\det\left(\begin{array}{ccccc}x_1 & -a_{1,2} & -a_{1,3}& \ldots & -a_{1,m} \\
-a_{2,1} & x_2 & -a_{2,3} & \ldots & -a_{2,m} \\
  & \vdots &   & \ddots & \vdots \\
-a_{m,1} & -a_{m,2} & -a_{m,3} & \ldots & x_m\end{array}\right).$$
Then for $x_i\ge d_i=\sum\limits_{j=1}^ma_{i,j},\,i =1,2,\ldots,m$
we have $D(x_1,x_2,\ldots,x_m)\ge0.$  The equality
$D(x_1,x_2,\ldots,x_m)=0$ holds if and only if $x_i=d_i,\,i
=1,2,\ldots,m.$
\end{lemma}

\begin{pf} We use induction on $m$ to prove the lemma. For $m=1$ we have $D(x_1)=x_1\ge a_{1,1}=0$ and $D(x_1)=0$ iff $x_1=a_{1,1}.$ For $m=2$ one has $D(x_1,x_2)=x_1x_2-a_{1,2}a_{2,1}\ge0$ with $D(x_1,x_2)=0$ if and only if $x_1=a_{1,2}$ and $x_2=a_{2,1}.$  Suppose that $m>2$ and lemma is true for all $D(x_1,x_2,\ldots,x_k)$ with $k<m.$

Denote by $D(x_1,\ldots,\hat{x}_i,\ldots,x_m),$ where $\hat{x}_i$
means that the variable $x_i$ is dropped, the $(i,i)$-th minor of
the matrix in the statement of lemma. We note that
$D^{\prime}_{x_i}(x_1,x_2,\ldots,x_k)=D(x_1,\ldots,\hat{x}_i,$ $\ldots,x_m)$, where $D^{\prime}_{x_i}(x_1,x_2,\ldots,x_k)$ is the partial derivative of $D(x_1,x_2,\ldots,x_k)$ with respect to the variable $x_i$.  Since
$$x_2\ge d_2=\sum\limits_{j=1}^ma_{2,j}\ge\sum\limits_{j=2}^ma_{2,j},x_3\ge d_3=\sum\limits_{j=1}^ma_{3,j}\ge\sum\limits_{j=2}^ma_{3,j}, \ldots,x_m\ge d_m=\sum\limits_{j=1}^ma_{m,j}\ge\sum\limits_{j=2}^ma_{m,j},$$ the function $D(x_2,\ldots,x_m)$ satisfies the conditions of lemma. Hence, $D(x_2,\ldots,x_m)\ge 0.$ In a similar way, for $i=2,\ldots,m$ we have
$$D^{\prime}_{x_i}(x_1,x_2,\ldots,x_m)=D(x_1,\ldots,\hat{x}_i\ldots,x_m)\ge0.$$ Since $D(d_1,d_2,\ldots,d_m)=0,$ we obtain $D(x_1,x_2,\ldots,x_m)\ge0$ for all $x_i\ge d_i,\,i =1,2,\ldots,m.$ If for some $i_0$ we have $x_{i_0}>d_{i_0},$ then, by induction, for all $i\neq i_{0}$ we get $D^{\prime}_{x_{i_0}}(x_1,x_2,\ldots,x_m)=D(x_2,\ldots,\hat{x}_{i_0}\ldots,x_m)>0$ and $D(x_1,x_2,\ldots,x_m)>0.$
\end{pf}

\bigskip
\begin{lemma}\label{lemma1}
We have  $P(z)=P(\frac{1}{z})$  and $$P(1)=0,\,P^{\prime}(1)=0,\,P^{\prime\prime}(1)< 0^{}.$$  In particular,  $P(z)$ has a double root $z=1.$
\end{lemma}

\begin{pf}
Since $P(1) =\det L(H) =0$. By~(\ref{newform}) we have $$P(z)=c_0
+\sum_{k=1}^{\ell(H)}\sum_{B \in {\cal C}_k(H)}c_B\left(\prod_{C \in
{\cal C}(B)}(z^{\alpha(C^+)} +z^{\alpha(C^-)})\right),$$ where
${\displaystyle c_B = (-1)^k\sum_{M \in {\cal
M}(H-B)}(-1)^{|E(M)|}\prod_{v \in V(H-B-M)}d(v)}$. Now we have
$$P^{\prime}(z)= \sum_{k=1}^{\ell(H)}\sum_{B \in {\cal C}_k(H)}c_B \left(\sum_{C \in
{\cal C}(B)} (\alpha(C^+)z^{\alpha(C^+)-1}
+\alpha(C^-)z^{\alpha(C^-)-1}) \prod_{C^{\prime} \in {\cal
C}(B-C)}(z^{\alpha(C^{\prime +})} +z^{\alpha(C^{\prime
-})})\right),$$ and hence $P^{\prime}(1)=0$ because  $\alpha(C^-) =
-\alpha(C^+)$. For any $C \in {\cal C}(H)$, let $L_{\bar{C}}$ be the
matrix obtained by deleting rows and columns of $L(H)$ corresponding
to vertices of $C$.  If $C$  is a Hamiltonian cycle, the matrix $L_{\bar{C}}$  is  empty. In this case  we set  $\det L_{\bar{C}}=1.$    By Lemma \ref{positive-det}, $\det L_{\bar{C}} >
0$. Furthermore, we have
\begin{eqnarray*}
\det L_{\bar{C}} &=& \sum_{B \in {\cal
B}(H-C)}(-1)^{|E(B_1)|}(-2)^{|{\cal C}(B)|}\prod_{v \in
V(H-C)\setminus V(B)}d(v)  \\
&=&  -\sum_{k=1}^{\ell(H)}2^{k-1}\sum_{B \in {\cal C}_k(H,C)}c_B.
\end{eqnarray*}
This implies that ${\displaystyle \sum_{k=1}^{\ell(H)}2^{k-1}\sum_{B
\in {\cal C}_k(H,C)}c_B} <0$. Note that

\begin{eqnarray*}
P^{\prime\prime}(z)&=& \sum_{C \in {\cal C}(H)} c_{C}
\left(\alpha(C^+)(\alpha(C^+)-1)z^{\alpha(C^+)-2}
+\alpha(C^-)(\alpha(C^-)-1)z^{\alpha(C^-)-2} \right) \\
&+& \sum_{k=2}^{\ell(H)}\sum_{B \in {\cal C}_k(H)}c_B
\left(\sum_{C\in {\cal C}(B)}\left(\alpha(C^+)(\alpha(C^+)-1)z^{\alpha(C^+)-2} \right. \right. \\
& & + \left. \alpha(C^-)(\alpha(C^-)-1)z^{\alpha(C^-)-2} \right) \times
\prod_{C^{\prime}
\in {\cal C}(B-C)}(z^{\alpha(C^{\prime +})} +z^{\alpha(C^{\prime -})}) \\
&+&  \sum_{C,C^{\prime} \in {\cal C}(B),C \neq C^{\prime}}
2\left((\alpha(C^+)z^{\alpha(C^+)-1}+\alpha(C^-)z^{\alpha(C^-)-1})
\times \right.\\
& &\left.\left(\alpha(C^{\prime +})z^{\alpha(C^{\prime+})-1}+
\alpha(C^{\prime -})z^{\alpha(C^{\prime -})-1}) \right)
\times \prod_{C^{\prime\prime} \in {\cal C}(B-C-C^{\prime})}(z^{\alpha(C^{\prime\prime +} )}
+z^{\alpha(C^{\prime\prime -})}) \right).
\end{eqnarray*}

So we have
\begin{eqnarray*}
P^{\prime\prime}(1)&=& \sum_{C \in {\cal C}(H)} 2 c_{C}
\alpha(C^+)^2  + \sum_{k=2}^{\ell(H)}\sum_{B \in {\cal C}_k(H)} 2^k
c_B   \sum_{C \in {\cal C}(B)} \alpha(C^+)^2 \\
&=& 2 \sum_{C \in {\cal C}(H)} \alpha(C^+)^2
\sum_{k=1}^{\ell(H)}2^{k-1}\sum_{B \in {\cal C}_k(H,C)}c_B =- 2
\sum_{C \in {\cal C}(H)} \alpha(C^+)^2 \det L_{\bar{C}}.
\end{eqnarray*}
Hence\begin{equation}\label{secondprime}P^{\prime\prime}(1) =  - 2
\sum_{C \in {\cal C}(H)} \alpha(C^+)^2 \det L_{\bar{C}}<0 \end{equation} and
\begin{equation}\label{secondprimemod}|P^{\prime\prime}(1)| =  2 {\displaystyle \sum_{C \in
{\cal C}(H)} \alpha(C^+)^2 \det L_{\bar{C}}}.\end{equation}

\end{pf}

\section{Main results}

The main results of the paper are the two following theorems.
\bigskip

\begin{theorem}\label{theorem1}
The number of spanning trees $\tau(n)$ in the graph $H_{n}$ is given by the formula
  $$\tau(n)= \frac{2n\,\tau(H)\varepsilon\,a_{0}^{n}}{P^{\prime\prime}(1)}\prod\limits_{\substack{P(z)=0\\z\neq1}}(z^n-1),$$
 where  the product is taken over all the roots different from $1$ of the voltage polynomial $P(z)$ of degree $2s,$\,$a_0$ is the  leading coefficient of  $P(z)$, $\varepsilon=(-1)^{s(n-1)}$ and $\tau(H)$ is the number of spanning trees in the graph $H.$\end{theorem}

\begin{pf} We note that $\tau(H_n)$  is non-negative integer. Then, by Proposition~\ref{theorem0}, we have
\begin{equation}\label{tauH}
\tau(H_n)=
\frac{\tau(H)}{n}\prod\limits_{j=1}^{n-1}P(\varepsilon_n^j).
\end{equation}
 Since $P(z)=P(\frac{1}{z}),$  Laurent polynomial has even degree, say $2s.$   Let $a_0$ be the leading coefficient of $P(z).$ To continue the proof we replace the Laurent polynomial $P(z)$ by a monic polynomial $\widetilde{P}(z)=a_0^{-1}z^{s}P(z).$
Then $\widetilde{P}(z)$ is a  polynomial of the degree $2s$ with the same roots as $P(z).$  Recall that by Lemma~~\ref{lemma1},  $P(z)$  has the double root $z=1.$

We note that
\begin{equation}\label{newP}
 \prod\limits_{j=1}^{n-1}P(\varepsilon_{n}^{j}) = \varepsilon\, a_{0}^{n-1} \prod\limits_{j=1}^{n-1}\widetilde{P}(\varepsilon_{n}^{j}),
\end{equation}
where $\varepsilon=(-1)^{s(n-1)}.$
By the basic properties of resultant we have
\begin{eqnarray}\label{Hlemma}
\nonumber &&\prod\limits_{j=1}^{n-1} \widetilde{P}(\varepsilon_{n}^{j})=\textrm{Res}(\widetilde{P}(z),\frac{z^{n}-1}{z-1}) =\textrm{Res}(\frac{z^{n}-1}{z-1},\widetilde{P}(z))\\
&&=\prod\limits_{z:\widetilde{P}(z)=0}\frac{z^{n}-1}{z-1}=\prod\limits_{z:{P}(z)=0}\frac{z^{n}-1}{z-1}= n^2\prod\limits_{{\substack{P(z)=0\\z\neq1}} }\frac{z^n-1}{z-1}.
\end{eqnarray}
  Combine (\ref{tauH}), (\ref{newP}) and (\ref{Hlemma}) we have the following formula for the number of spanning trees
\begin{equation}\label{modulo}
\tau(H_n)=n\,\tau(H)\varepsilon\, a_{0}^{n-1}\prod_{{\substack{P(z)=0\\z\neq1}} } \frac{z^n-1}{z-1}=
n\,\tau(H)\varepsilon \,a_{0}^{n-1}\prod\limits_{{\substack{P(z)=0\\z\neq1}} } (z^n-1)\big/\prod_{{\substack{P(z)=0\\z\neq1}} }(z-1).
\end{equation}

We note that
$P(z)=  z^{-s}(z-1)^2R(z),$ where $R(z),\,R(1)\neq0$ is an integer polynomial of even degree $2s-2$ with the leading coefficient $a_0.$ By the direct calculation  we obtain  $P^{\prime\prime}(1)=2R(1).$  Then

\begin{equation}\label{eqQs}\prod_{{\substack{P(z)=0\\z\neq1}} }(z-1)=\prod_{{\substack{R(z)=0}} }(z-1)=\frac{1}{a_0}R(1)=\frac{1}{2a_0}P^{\prime\prime}(1).
\end{equation}
Substituting equation (\ref{eqQs}) into equation (\ref{modulo}) we finish the proof of the theorem.

\bigskip
\end{pf}

One can halve the number of multiples in Theorem~\ref{theorem1} by making use of the following trick. By Lemma~\ref{lemma1}, we know that the Laurent polynomial $P(z)$ of degree $2s$ satisfies the property $P(z)=P(\frac{1}{z}).$ Hence, it can be written in the form
$$P(z)=a_0(z^s+\frac{1}{z^s})+a_{1}(z^{s-1}+\frac{1}{z^{s-1}})+\ldots+a_{s-1}(z+\frac{1}{z})+a_s.$$ We note the following identity $\frac12(z^n+\frac{1}{z^n})=T_n(\frac12(z+\frac{1}{z})),$  where $T_n(x)=\cos(n\arccos(x))$ is the Chebyshev polynomial of the first kind. Let $w=\frac12(z+\frac{1}{z}),$  then $P(z)=Q(w),$ where
\begin{equation}\label{polynomialQ}\
Q(w)=2a_0T_s(w)+2a_1T_{s-1}(w)+\ldots+2a_{s-1}T_{1}(w)+a_s.
\end{equation}
  We will call $Q(w)$ by a $\it Chebyshev \,\,transform$  of $P(z).$  Note that $Q(w)$  is an order $s$  polynomial with  the leading coefficient  $\widetilde{a}_0=2^sa_0.$

  By Lemma~\ref{lemma1}, the roots of polynomial $P(z)$  are $1,1,z_1,\frac{1}{z_1},\ldots,z_{s_1},\frac{1}{z_{s-1}},$ where $\,z_i\neq1.$ Then the roots of polynomial $Q(w)$  are given by the list $1,w_1,\ldots,w_{s-1},$ where $w_i=\frac12(z_i+\frac{1}{z_i})$  and $w_i\neq1.$  By direct calculation and Lemma~\ref{lemma1} we obtain $Q(1)=0,\,Q^{\prime}(1)=P^{\prime\prime}(1)\neq 0.$  Also we have

  $$(z_p^n-1)(z_p^{-n}-1)=2-(z_p^n+z_p^{-n})=2-2T_n(w_p),\,p=1,2,\ldots,s-1.$$
  We observe that  the product $\prod\limits_{\substack{P(z)=0\\z\neq1}}(z^n-1)$  in Theorem~\ref{theorem1} can be rewritten as

  \begin{equation}\label{rewrite}
  \prod\limits_{p=1}^{s-1}(z_p^n-1)(z_p^{-n}-1)= \prod\limits_{p=1}^{s-1}(2-2T_n(w_p)).
  \end{equation}

 This gives a way to represent Theorem~\ref{theorem1}  in the following form.

 \begin{theorem}\label{theorem2}Let $Q(w)$  be the $\it Chebyshev \,\,transform$  of the voltage polynomial $P(z)$ of a cyclic covering $H_n.$ Then the  number of spanning trees $\tau(n)$ in the graph $H_{n}$ is given by the formula
  $$\tau(n)=\frac{2 n\tau(H)\varepsilon\,a_{0}^{n}}{P^{\prime\prime}(1)}\prod\limits_{p=1}^{s-1}(2-2T_n(w_p)),$$
 where the product is taken over all the roots different from $1$ of the polynimial $Q(w)$  of degree $s$, $\varepsilon=(-1)^{s(n-1)}$, ${a}_0$  is the  leading coefficient of  $P(z).$   \end{theorem}

We note that   $P^{\prime\prime}(1)$ is given by the formula~(\ref{secondprime}).
 \section{Generating function for the number of spanning trees }

In this section, we investigate the properties of the generating function for the numbers of spanning trees of
  graph $H_{n}.$  Our aim is to prove the following result.
\bigskip

\begin{theorem}\label{theoremR1} Let $\tau(n)$ be the number of spanning trees in the   graph $H_{n}.$ Then $$F(x)=\sum\limits_{n=1}^\infty\tau(n)x^n$$ is a rational function with integer coefficients. Moreover, $F(a_0x)=F(\frac{1}{a_0 x}),$ where ${a}_0$ is the  leading coefficient of the voltage polynomial $P(z).$ The latter gives a way to represent $F(x)$ as a rational function of $u=\frac12(a_0x+\frac{1}{a_0 x}).$
\end{theorem}
\bigskip
The proof of   Theorem~\ref{theoremR1}  is based on the following proposition.
\begin{prop}\label{proposition1} Let $R(z)$  be a degree $2s$ polynomial with integer coefficients.  Suppose that all  the roots  of the polynomial  $R(z)$  are $\xi_1,\xi_2,\ldots,\xi_{2s-1},\xi_{2s}.$   Then $$F(x)=\sum\limits_{n=1}^\infty n\prod\limits_{j=1}^{2s}(\xi_{j}^{n}-1)x^n$$  is  a rational  function with integer coefficients.

Moreover, if  $\xi_{j+s}=\xi_{j}^{-1},\,j=1,2,\ldots,s,$  then $F(x)=F(1/x).$
\end{prop}

\begin{pf} First of all, we note that $F(x)=x\frac{d G(x)}{d x},$
where $$G(x)=\sum\limits_{n=1}^\infty\prod\limits_{j=1}^{2s}(\xi_{j}^{n}-1)x^n.$$

Denote by $\sigma_k=\sigma_k(x_1,x_2,\ldots,x_{2s})$ the $k$-th basic symmetric polynomial in variables $x_1,x_2,\ldots,x_{2s}.$  Namely,  $$\sigma_0=1,\,\sigma_1=x_1+x_2+\ldots+x_{2s}, \,\sigma_2=x_1x_2+x_1x_3+\ldots+x_{2s-1}x_{2s},\ldots, \sigma_{2s}=x_1x_2\ldots x_{2s}.$$  Then
$$G(x)=G_{2s}(x)-G_{2s-1}(x)+\ldots-G_1(x)+ G_0(x),$$  where $$G_k(x)= \sum\limits_{n=1}^\infty\sigma_k(\xi_1^n,\xi_2^n,\ldots, \xi_{2s}^n)x^n,\,k=0,1,\ldots,2s.$$  We have
$\sigma_k(\xi_1^n,\xi_2^n,\ldots, \xi_{2s}^n)=\sum\limits_{\substack{1\le{j_1}<{j_2}<\ldots<{j_k}\le2s}}\xi_{j_1}^n\xi_{j_2}^n\ldots\xi_{j_k}^n.$  Hence,

$$G_k(x)=\sum\limits_{\substack{1\le{j_1}<{j_2}<\ldots<{j_k}\le2s}}\frac{\xi_{j_1}\xi_{j_2}\ldots\xi_{j_k}}{1-\xi_{j_1}\xi_{j_2}\ldots\xi_{j_k} x}$$   and \begin{equation}\label{derivative}
F_k(x)=x\frac{d G_k(x)}{d x}=
\sum\limits_{\substack{1\le{j_1}<{j_2}<\ldots<{j_k}\le2s}}\frac{ \xi_{j_1}\xi_{j_2}\ldots\xi_{j_k} x}{(1-\xi_{j_1}\xi_{j_2}\ldots\xi_{j_k} x)^2}.
\end{equation}  We note that $F_k(x)$ is a   symmetric function in the roots $\xi_1,\xi_2,\ldots,\xi_{2s}$  of the integer  polynomial $R(x).$    By the Vieta theorem, $F_k(x)$ is   a rational function with integer coefficients. Since

\begin{equation}\label{alternative} F(x)=F_{2s}(x)-F_{2s-1}(x)+\ldots-F_1(x)+ F_0(x),\end{equation}
the same is true for $F(x).$

To prove the second statement of the proposition,  consider the product $\xi=\xi_{j_1}\xi_{j_2}\ldots\xi_{j_k}.$
Since $\xi_{j+s}=\xi_{j}^{-1},$   the term $\frac{\xi x}{(1-\xi x)^2}$ comes into  the sum   (\ref{derivative}) together with the term
$\frac{  \xi^{-1}x }{(1-\xi^{-1}x)^2}.$ One can check   that  the function  $\varphi(x)=\frac{\xi x}{(1-\xi x)^2}+\frac{\xi^{-1}x}{(1-\xi^{-1}x)^2} $ satisfies the condition $\varphi(x)=\varphi(1/x).$  Then, for   $k=0,1,\ldots, 2s,$  we have $F_k(x)=F_k(1/x).$   Finally, by (\ref{alternative}) we obtain $F(x)=F(1/x).$

\end{pf}

{\bf Proof of Theorem~\ref{theoremR1}.}
\medskip\noindent  We employ Theorem~\ref{theorem1}  and Proposition~\ref{proposition1} to prove the theorem.
 Consider the polynomial $R(z)=z^{s}P(z)/(z-1)^2,$  where $P(z)$  is the voltage polynomial of the graph  $H_n.$  Note that $R(z)$ is a polynomial of order
 $2s-2$ with integer coefficients.
 Recall that all the roots of the polynomial $R(z)$ are the roots of  $P(z)$ different from $1.$
  We will use Theorem~\ref{theorem1} to find $\tau(n)$.  Then
\begin{equation}\label{gool}\tau(n)=\frac{2n\tau(H)\varepsilon\,a_0^n}{P^{\prime\prime}(1)}\prod\limits_{\substack{P(z)=0\\z\neq1}}(z^n-1)=
\frac{2n\tau(H)\varepsilon\,a_0^n}{P^{\prime\prime}(1)}\prod_{p=1}^{s-1}(z_p^n-1)(z_p^{-n}-1),
\end{equation}
  where $z_p,\,z_p^{-1},\, p=1,2,\ldots,s-1$
 are all the roots  of the   polynomial $R(z)$  and $\varepsilon=(-1)^{s(n-1)}.$

Then, by virtue of (\ref{gool}),  the generating function $F(x)=\sum_{n=1}^\infty\tau(n)x^n$  can be represented in the form
$$F(x)=(-1)^{s}\frac{2\tau(H)}{P^{\prime\prime}(1)}\sum\limits_{n=1}^\infty n\prod_{p=1}^{s-1}(z_p^n-1)(z_p^{-n}-1)((-1)^sa_0x)^n,$$
Since  $a_0$  and $\frac{2\tau(H)}{P^{\prime\prime}(1)}$  are rational  numbers, by Proposition~\ref{proposition1}, $F(x)$  is a rational function with integer coefficients satisfying $F((-1)^sa_0x)=F(\frac{1}{(-1)^sa_0x}).$ Hence, $F(a_0x)=F(\frac{1}{a_0x}).$ \hfill\qed

\section{Asymptotic formulas for the number of spanning trees }

In this section we aim to find an asymptotic formula for $\tau(H_n)$.
The following theorem is the main result of this section.

\bigskip

\begin{theorem}\label{theorem3}
The asymptotic behavior  for number of spanning trees in the graph $H_n $   is given by the formula
$$\tau(H_n)\sim \frac{2n\,\tau(H)}{q}A^n,\,n\to\infty,$$ where $q=|P^{\prime\prime}(1)|$  and
$A=\exp({\int\limits_{0}^{1}\log|P(e^{2 \pi \mathtt{i} t})|\textrm{d}t})$  is the Mahler measure of the Laurent polynomial $P(z).$
\end{theorem}

Note that formula~(\ref{secondprimemod}) allows us to find the value $|P^{\prime\prime}(1)|$ explicitly. In order to prove the theorem, preliminary, we establish the following lemma.
\bigskip
\begin{lemma}\label{nounityroots}
For any complex number $e^{\mathtt{i} \varphi}$ with $\varphi\in\mathbb{R}$, one has $P(e^{\mathtt{i} \varphi})\ge0.$  Furthermore,  $P(e^{\mathtt{i} \varphi})=0$ if and
only if $e^{\mathtt{i} \varphi}=1$.  \end{lemma}

\begin{pf}
Let   $L(z)$ be the matrix
\begin{equation}
\left(\begin{array}{cccc}
p_{11}(z) & p_{12}(z)  & \ldots & p_{1r}(z) \\
p_{21}(z) & p_{22}(z)  & \ldots & p_{2r}(z) \\
  \vdots &   & \ddots & \vdots \\
p_{r1}(z) & p_{r2}(z)  & \ldots & p_{rr}(z)\\
\end{array}\right),\end{equation} where $p_{ii}(z)=a_{i,i}+\sum\limits_{j=1}^ra_{i,j}-\sum\limits_{k=1}^{a_{i,i}}(z^{\alpha_{k}(i,i)}  + z^{-\alpha_{k}(i,i)} )$  and
$p_{ij}(z)=-\sum\limits_{k=1}^{a_{i, j}}z^{\alpha_k(i,j)}$ for $i
\neq  j.$ Now $P(z)=\det L(z)$. For any complex number $e^{\mathtt{i} \varphi}$, $L(e^{\mathtt{i} \varphi})$ is a Hermitian matrix and so it suffices to
show that for any complex vector ${\mathbf x} \in \mathbb{C}^r$, ${\mathbf
x}^H L(e^{\mathtt{i} \varphi}){\mathbf x}  \ge 0$ and furthermore, ${\mathbf
x}^H L(e^{\mathtt{i} \varphi}){\mathbf x}  = 0$ for some ${\mathbf x} \in \mathbb{C}^r
\setminus \{{\mathbf 0}\}$ if and only if $e^{\mathtt{i} \varphi}=1$. Let
$D(z)$ be an $r \times r$ diagonal matrix whose $(i,i)$-entry is
$d_{ii}(z) = 2a_{i,i}
-\sum\limits_{k=1}^{a_{i,i}}(z^{\alpha_{k}(i,i)}  +
z^{-\alpha_{k}(i,i)} )$. Set $L_1 (z) = L(z)-D(z)$. Note that $D(1)$
is a zero matrix and $L_1(1)$ is the Laplacian matrix of the graph
obtained by deleting all loops of $H$. So $P(1) = \det L(1) =0$.

Let $v_1, \ldots, v_r$ be all vertices of $H$ and for any $k=1,
\ldots, a_{i, j}$ with $i \neq j$, let $e_k\{i,j \}$ be the edge
between $v_i$ and $v_j$ such that the voltages assigned to two
directed edges induced by $e_k\{i,j \}$ are $\alpha_k(i,j)$ and
$-\alpha_k(i,j)$. Let $m$ be the number of edges of $H$ which is not
a loop. Denote by $I(z)$  the $m \times r$ matrix whose rows and columns
are labeled by non-loop edges and vertices of $H$, respectively,
such that the entry of $I(z)$ corresponding to the row $e_k\{i,j \}$
and column $v_t$ is $z^{-\frac{\alpha_k(i,j)}{2}}$ if $t=i$;
$-z^{\frac{\alpha_k(i,j)}{2}}$ if $t=j$; and $0$ otherwise. Now  for
any complex number $e^{\mathtt{i} \varphi}$, $L_1(e^{\mathtt{i} \varphi}) =
I(e^{\mathtt{i} \varphi})^H I(e^{\mathtt{i} \varphi})$. On the other hand,
$D(e^{\mathtt{i} \varphi})$ is the diagonal matrix such that $d_{ii}(e^{\mathtt{i} \varphi
}) = 2\sum\limits_{k=1}^{a_{i,i}}(1-\cos(\alpha_{k}(i,i)\varphi))$.
Therefore, for any complex vector ${\mathbf x} \in \mathbb{C}^r$ with
${\mathbf x}^H = (\bar{x}_1, \ldots, \bar{x}_r)$,
\begin{eqnarray*}
{\mathbf x}^H L(e^{\mathtt{i} \varphi}){\mathbf x} &=& {\mathbf x}^H
D(e^{\mathtt{i} \varphi}){\mathbf x}  + {\mathbf x}^H L_1(e^{\mathtt{i} \varphi}){\mathbf x} \\
&=& 2\sum_{k=1}^{a_{i,i}}(1-\cos(\alpha_{k}(i,i)\varphi))|x_i|^2
+2\sum_{i,j, i \neq j} \sum_{k=1}^{a_{i, j}} \left|
x_ie^{-\frac{\alpha_k(i,j)\mathtt{i}\varphi}{2}}
-x_je^{\frac{\alpha_k(i,j)\mathtt{i}\varphi}{2}} \right|^2 \ge 0.
\end{eqnarray*}
This implies that $P(e^{\mathtt{i} \varphi}) = \det L(e^{\mathtt{i} \varphi}) \ge 0$.
Furthermore, ${\mathbf x}^H L(e^{\mathtt{i} \varphi}){\mathbf x}=0$ if and
only if $(1-\cos(\alpha_{k}(i,i)\varphi))|x_i|^2 =0$ for any
$k=1,\ldots, a_{i,i}$ and $x_i= x_je^{\alpha_k(i,j)\mathtt{i} \varphi}$ for
any $k =1, \ldots, a(i,j)$ with $i \neq j$. Note that if $x_i =0$
for some $i=1,\ldots,r$, then $x_j=0$ for all $j=1, \ldots, r$,
namely, ${\mathbf x} ={\mathbf 0}$ because  $H$ is connected.

Assume that $e^{\mathtt{i} \varphi} \neq 1$ and ${\mathbf x}^H L(e^{\mathtt{i} \varphi}){\mathbf x}=0$.  Suppose that ${\mathbf x} \neq
{\mathbf 0}$, and equivalently  $x_j\neq 0$ for any $j=1, \ldots,
r$. If there is a loop incident to $v_i$ such that $\alpha_{k}(i,i)
\neq 0$ for some $k$, then $x_i =0$, a contradiction. So we may
assume that for any loop in $H$, its voltage is $0$, and hence
$D(z)$ is a zero matrix. Now for any edge $e_k\{i,j\}$, we have $x_i
= x_je^{\alpha_k(i,j)\mathtt{i}\varphi}$. This implies that for any cycle
$C$ containing $v_i$, we have $x_i = x_ie^{\alpha(C^+)\mathtt{i}\varphi}$.
Since $x_i \neq 0$, we have $e^{\alpha(C^+)\mathtt{i}\varphi}=1$. Hence for
any cycle $C$ in $H$, $e^{\alpha(C^+)\mathtt{i} \varphi}=1$.   Since there exist cycles $C_1,  \ldots, C_s$ such that
$k_1\alpha(C_1^+)+\cdots +k_s\alpha(C_s^+)=1$ for some integers $k_1,\ldots, k_s$, we have
$$1 = e^{k_1\alpha(C_1^+)\mathtt{i}\varphi} \cdots e^{k_s\alpha(C_s^+)\mathtt{i} \varphi
} = e^{(k_1\alpha(C_1^+)+\cdots +k_s\alpha(C_s^+))\mathtt{i}\varphi }
=e^{\mathtt{i} \varphi},$$ which is a contradiction. So for any ${\mathbf x}
\in \mathbb{C}^r \setminus \{{\mathbf 0}\}$, ${\mathbf x}^H L(e^{\mathtt{i}\varphi
}){\mathbf x}  > 0$, and hence $P(e^{\mathtt{i} \varphi}) = \det
L(e^{\mathtt{i} \varphi}) >0$.
\end{pf}

Now we are  able to prove Theorem~\ref{theorem3}.

\begin{pf}
 By Lemmata~\ref{lemma1}   and \ref{nounityroots},   the polynomial $P(z)$ has the roots  $1,1,z_{j}$ and $1/z_{j},\,j=1,2,\ldots,s-1$ with the property $|z_{j}|\neq1,\,j=1,2,\ldots,s-1.$ By theorem \ref{theorem3} we have $$\tau(n)=\frac{2 n\tau(H)|a_0|^n}{q}\prod\limits_{j=1}^{s-1}|(z_j^n-1)(z_j^{-n}-1)|.$$  Replacing $z_{j}$ by $1/z_{j},$ if it is necessary, we can assume that  $|z_j|>1$ for all $j=1,2,\ldots,s-1.$ Then  and $|z_j^n-1|\sim|z_{s}|^{n}$ and $|z_j^{-n}-1|\sim1$ as $n\to\infty.$ Hence
$$\frac{2\,n \tau(H)|a_0|^n}{q}\prod_{j=1}^{s-1}|(z_j^n-1)(z_j^{-n}-1)|\sim\frac{2\,n \tau(H)|a_0|^n}{q}\prod_{j=1}^{s-1} |z_{j}|^{n}=\frac{2n\,\tau(H)}{q}A^n,$$ where
 $A=|a_0|\prod_{j=1}^{s-1} |z_{j}|^{n}=|a_0|\prod\limits_{ P(z)=0,\,|z|>1 }|z|$ is the Mahler measure of the polynomial $P(z).$ By (\cite{EverWard}, p.~67), we have $A=\exp\left(\int\limits_{0}^{1}\log|P(e^{2 \pi \mathtt{i} t })|\textrm{d}t\right)$.
\end{pf}

As a consequence  of Theorem~\ref{theorem3}, we obtain the following theorem.

\begin{theorem}\label{thermolimit} The thermodynamic   limit of the sequence of graphs $H_{n}$  is equal to the small Mahler  measure of the Laurent polynomial $P(z).$  More precisely,

$$\lim\limits_{n\to\infty}\frac{\log \tau(H_{n})}{n}= \int\limits_{0}^{1}\log|P(e^{2\pi\mathtt{i} t})|d t.$$
\end{theorem}

 \section{Examples}

\begin{enumerate}
\item{\textbf{Haar graph $H_n(1,2)$}.}
($n$-fold cyclic covering of the theta graph with voltages $0,1,2,$ see Fig. 1.) We get  $P(z)=-(1-z)^2(1+4z+4z^2)/z^2$  and $Q(w)=4(1-w)(2+w).$ By Theorem~\ref{theorem2}, we have   $\tau(n)=(-1)^{n}(T_n(-2)-1).$
By Theorem~\ref{theorem3}, $\tau(n)$ has the following asymptotic
$\tau(n)\sim\frac{n}{2}(2+\sqrt{3})^n,\, n\to\infty.$ Also,
$$\sum\limits_{n=1}^{\infty}\tau(n)x^n=\frac{3 (1 + u + u^2)}{2 (1 - u)(2 + u)^2 },$$ where $u=-\frac12(x+1/x).$

 \begin{figure}[h]
\begin{center}
  \includegraphics[scale=1.55]{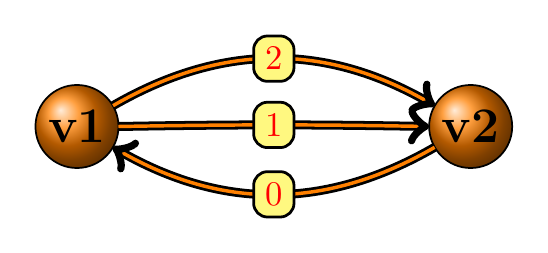}

 $$\text{Voltage scheme  of the Haar graph  }  H_n(0,1,2) $${}$$Figure \,\,1.$$
\end{center}

\end{figure}

\item{\textbf{Cyclic coverings of the three circles graph}.} (See Fig. 2  for the voltage scheme.) We have $$P(z)=14 - 4 (z^{\alpha} +
 z^{-\alpha}) - 4 (z^{\beta} +
 z^{-\beta})-   (z^{\gamma} +
 z^{-\gamma})  +( z^{\alpha+\beta} +z^{-\alpha- \beta}) +( z^{\alpha-\beta} +z^{-\alpha+\beta}).$$ In this case, $P^{\prime\prime}(1)=-2(2\alpha^2+2\beta^2+ \gamma^2)$  and Chebyshev transform of  $P(w)$  has the form
 $$Q(w)=14-8T_{\alpha}(w)-8T_{\beta}(w)-2T_{\gamma}(w)+2T_{\alpha+\beta}(w)+2T_{\alpha-\beta}(w).$$
  In particular, if $(\alpha,\beta,\gamma)=(1,2,3),$
 then $Q(w)=2(1-w)(11+8w).$  Here, the number of spanning trees in the graph $H_n$ is given by the formula
 $$\tau(n)=\frac{n(-4)^{n+1}}{19}(1-T_n(-\frac{11}{8})).$$  Also, we have
 $$\sum\limits_{n=1}^\infty\tau(n)x^n=\frac{2(3+8u+8u^2)}{(1-u)(11+8u)^2},$$  where $u=-\frac12(4x+\frac{1}{4x}).$

 \begin{figure}[h]
\begin{center}
  \includegraphics[scale=1.0]{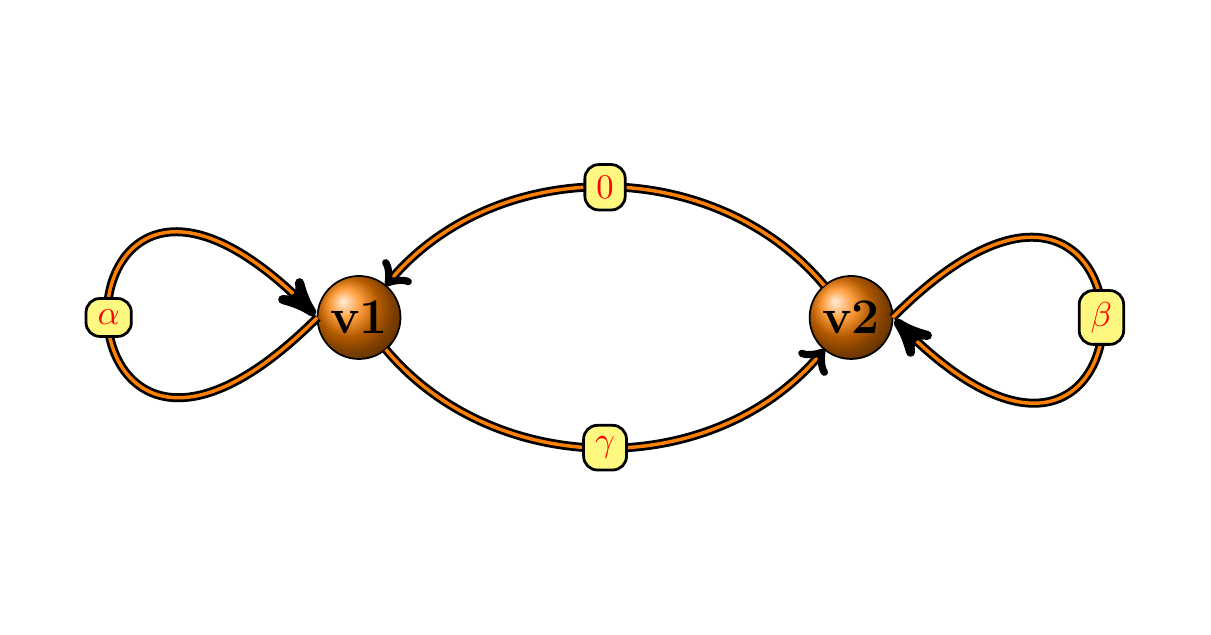}
 \vspace{-5mm}
 $$\text{Voltage scheme of the three circles graph }$${}$$Figure\,\, 2.$$
\end{center}

\end{figure}

 \item{\textbf{Cyclic coverings of the triangle with a loop and a double edge}.}  Let $H$ be the graph on vertices $v_1,v_2,v_3$ with edges $v_1v_2,\,v_1v_3,\,v_2v_3,\,v_2v_3$ and the loop $v_1v_1.$ We attach to this graph the following voltages
$\alpha_1(1,1)=\alpha,\,\alpha_1(1,2)=0,\,\alpha_1(1,3)=0,
\,\alpha_1(2,3)=\beta,\,\alpha_2(2,3)=-\gamma.$ (See Fig. 3.) We get the following voltage polynomial

\begin{figure}[h]
\begin{center}
  \includegraphics[scale=1.0]{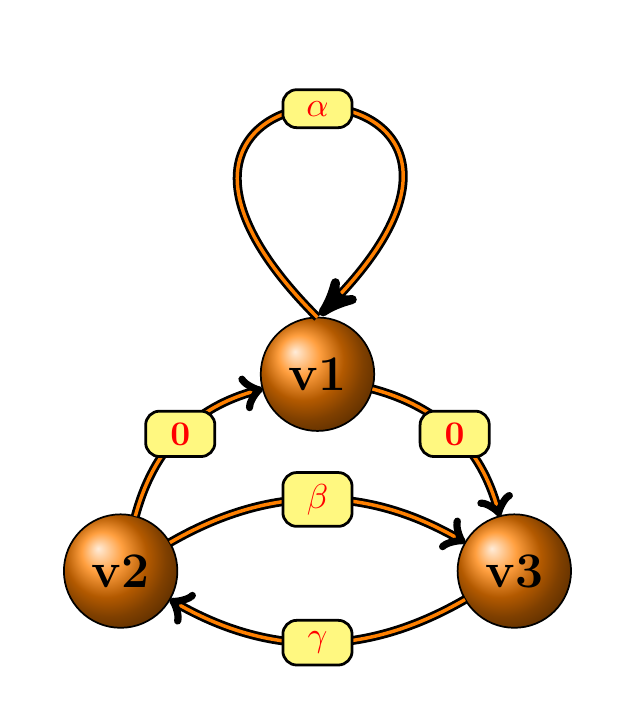}

 $$\text{The triangle with a loop and a double edge }$${}$$Figure\,\, 3.$$
\end{center}

\end{figure}

$P(z)=22-7(z^{\alpha}+z^{-\alpha})-(z^{\beta}+z^{-\beta})-(z^{\gamma}+z^{-\gamma})
-4(z^{\beta+\gamma}+z^{-\beta-\gamma})+(z^{\alpha-\beta-\gamma}+z^{-\alpha+\beta-\gamma}) + (z^{\alpha+\beta+\gamma}+z^{-\alpha-\beta-\gamma}).$

In this case, $P^{\prime\prime}(1)=-2(5\alpha^2+\beta^2+\gamma^2+2(\beta+\gamma)^2).$
The Chebyshev transform of $P(z)$  is given by  the formula

$$Q(w)=22-14T_{\alpha}(w)-2T_{\beta}(w)-2T_{\gamma}(w)-8T_{\beta+\gamma}(w)+2T_{\alpha-\beta-\gamma}(w)+2T_{\alpha+\beta+\gamma}(w).$$In particular, if $(\alpha,\beta,\gamma)=(1,2,-3),$
 then $Q(w)=8(1-w)(3+w+w^2).$  Now, $$\tau(n)=\frac{3n}{5}(1- T_n(\frac{-1-i\sqrt{11}}{2}))(1-T_n(\frac{-1+i\sqrt{11}}{2})).$$

\item{\textbf{{Y}-graphs $Y(n;k,l,m)$}.} The notion of an $Y$-graph was introduced in \cite{Biggs} and \cite{BoChMo}. Its symmetry properties were investigated in \cite {HorBou}.

 The voltage scheme of the graph $Y(n;k,l,m)$ is shown on Fig. 4. In this case we have
$Q(w)=3A B C-A B-B C-A C,$ where
$A=3-2T_k(w),\,B=3-2T_l(w),\,C=3-2T_m(w).$  In particular, for the graph $Y(n;1,1,1)$ one has $Q(w)=6(1-w)(3-2w)^2.$
By Theorem~\ref{theorem2}, the number $\tau(n)$ of spanning trees in
the   graph $Y(n;1,1,1)$ is given by the formula
$$\tau(n)=n3^{n-1}(2-2T_n(\frac32))^2.$$ Hence,
 $\tau(n)=3^{n-1}n \,L_n^4, $ if   $n$ is odd, and
  $\tau(n)=25\cdot3^{n-1}n\, F_n^4,$ if   $n$ is even,
 where  $L_n$  and $F_n$  are the Lucas  and the Fibonacci numbers respectively.  Also,
 $$\sum\limits_{n=1}^\infty\tau(n)x^n=\frac{-575-322u+298u^2+112u^3-8u^4}{3(1+u)(21+20u+4u^2)^2},$$  where $u=-\frac12(3x+\frac{1}{3x}).$

 \begin{figure}[h]
\begin{center}
  \includegraphics[scale=0.5]{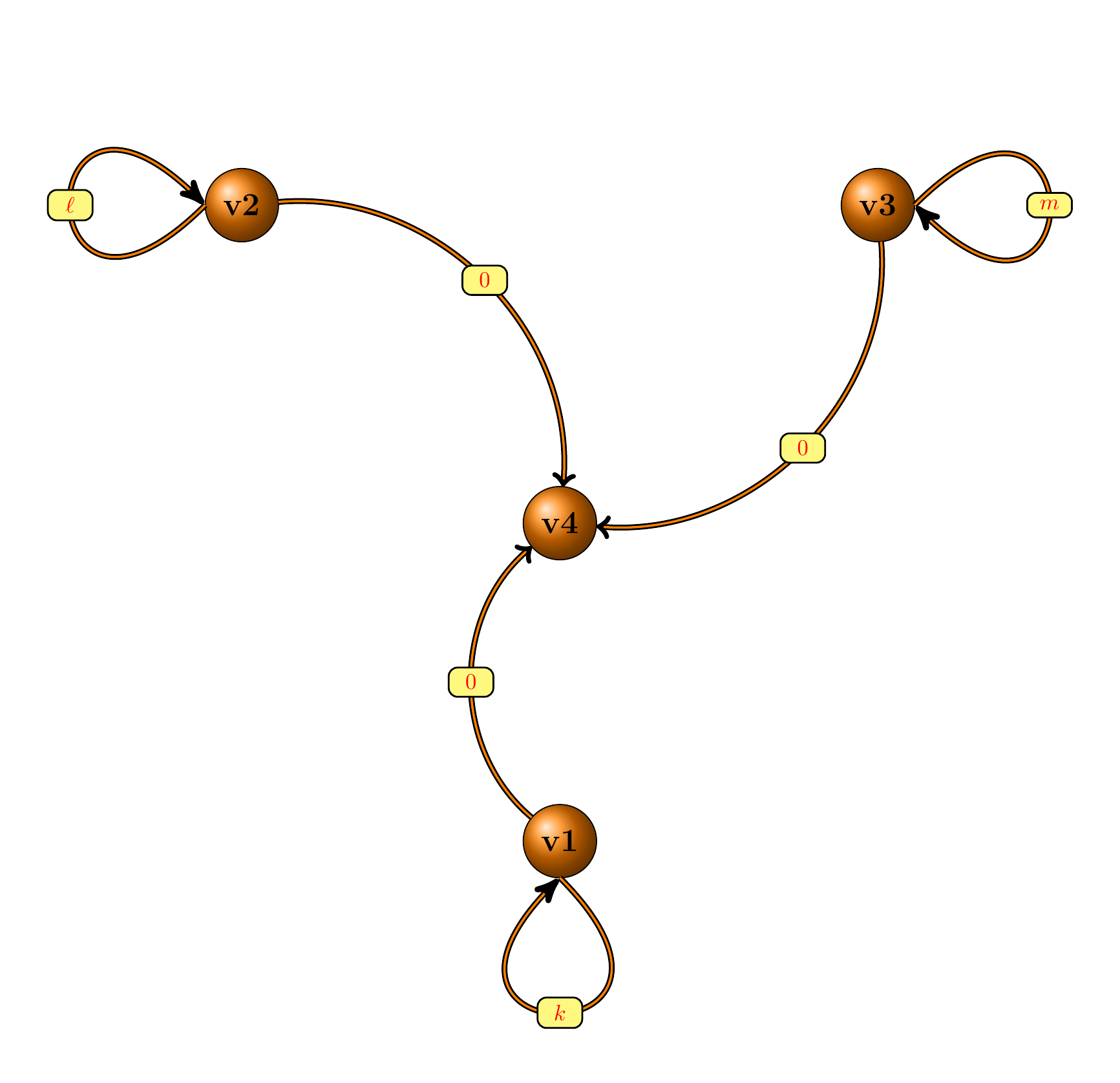}

 $$\text{Voltage scheme of the   graph }  Y(n;k,l,m) $${}$$Figure\,\, 4.$$
\end{center}

\end{figure}

\end{enumerate}

\section*{ACKNOWLEDGMENTS}
 The author are thankful to Professor Sergei  Lando who suggested us to prove Theorem~\ref{theoremR1}. The work was partially supported by the Korean-Russian bilateral project.  The first author was supported in part by the Basic Science Research
Program through the National Research Foundation of Korea (NRF)
funded by the Ministry of Education (2018R1D1A1B05048450).  The  second and the third  authors
were are partially supported by the Russian Foundation for Basic Research (grants 18-01-00420 and 18-501-51021).  The results given
in Sections 3 and 4 are supported by the Laboratory of Topology and Dynamics, Novosibirsk State University
(contract no. 14.Y26.31.0025 with the Ministry of Education and Science of the Russian Federation).

\end{document}